\documentclass[11pt,leqno]{article}
\normalfont
\renewcommand{\baselinestretch}{1.15}
\usepackage{amsfonts}

\newfont{\eulercursive}{eurm10 at 11pt}

\newcommand{\QED}{\raisebox{0.5mm}{\fbox{\rule{0mm}{1.5mm}\ }}}

\newcounter{myfn}[page]

\newcommand{\NotMarsFriendlyLemma}{Lemma 2.5}

\newcommand{\ThreeCycleLemma}{Proposition 1}

\begin{document}

\newpage
\setcounter{page}{1} 
\renewcommand{\baselinestretch}{1}

\vspace*{-0.45in}

\begin{center}
{\large \bf Eriksson's numbers game on certain edge-weighted 
three-node cyclic graphs}

Robert G.\ Donnelly

Department of Mathematics and Statistics, Murray State
University, Murray, KY 42071
\end{center}

\begin{abstract}
The numbers game is a one-player game played on  
a finite simple 
graph with certain ``amplitudes'' assigned to its edges 
and with an initial assignment of real 
numbers to its nodes.  The moves of the game successively transform 
the numbers at the nodes using the amplitudes in a certain way.  
This game and its interactions with Coxeter/Weyl group theory 
and Lie theory  
have been studied by many authors.  
Following Eriksson, we  
allow the amplitudes on graph edges 
to be certain real numbers. Games played on such graphs are 
``E-games.''  
We show that for certain such three-node cyclic graphs, any numbers game 
will diverge when 
played from an initial assignment of nonnegative real numbers not all 
zero. This result is a key step in a Dynkin diagram classification 
(obtained elsewhere) of all 
E-game graphs which meet a certain finiteness requirement. 
\begin{center}

\ 
\vspace*{-0.1in}


\end{center}
\end{abstract}

\vspace{1ex} 
The numbers game is a one-player game played on a finite simple graph 
with weights (which we call ``amplitudes'') on its edges 
and with an initial assignment of 
real numbers 
to its nodes.  
Each of the two edge amplitudes (one for each direction) 
will be certain negative real numbers. 
The 
move a player can make 
is to ``fire'' one of the nodes with a positive number.  This move 
transforms the number at the fired node 
by changing its sign, and it also 
transforms the number at each adjacent node in a certain way 
using an amplitude   
along the incident edge.  
The player fires the nodes in some sequence of 
the player's choosing, continuing until no node has a positive 
number.  
This game was formulated by Mozes \cite{Mozes} for graphs 
with integer amplitudes.  
It has also been studied by Proctor 
\cite{PrEur}, \cite{PrDComplete}, Bj\"{o}rner 
\cite{Bjorner}, 
and Wildberger 
\cite{WildbergerAdv}, \cite{WildbergerEur}, \cite{WildbergerPreprint}.  
The game is the subject of \S 4.3 of the book \cite{BB} by 
Bj\"{o}rner and Brenti. 
The numbers game 
facilitates computations with Coxeter groups and their geometric 
representations (e.g.\ see \S 4.3 of \cite{BB}).  
See \cite{DonEnumbers} for discussion of further connections and 
applications, where  
mainly we draw on Eriksson's ground-breaking work in  \cite{ErikssonThesis}, 
\cite{ErikssonDiscrete}, and \cite{ErikssonEur}. 

Our purpose here is to show in \ThreeCycleLemma\ that  
certain three-node cyclic ``E-GCM graphs'' (see \cite{DonEnumbers} 
for a definition) are not ``admissible'': that is, 
any numbers game played on such a graph from 
a nontrivial initial assignment of nonnegative numbers will not  
terminate.  
Our main interest in this proposition is that it furnishes 
a key step for the proof given in \cite{DonEnumbers} of the following Dynkin diagram 
classification result: 
A connected E-GCM graph  
has a nontrivial initial assignment of nonnegative 
numbers such 
that the numbers game terminates in a finite number of steps 
if and only if it is a connected ``E-Coxeter graph'' 
corresponding to an irreducible finite Coxeter group.  
(Another proof of this classification result is given in \cite{DE}.)  
All further motivation, 
definitions, and 
preliminary results needed to understand the statement and  
proof of \ThreeCycleLemma\ are given 
in \cite{DonEnumbers}.

\noindent 
{\bf \ThreeCycleLemma}\ \ 
{\sl Suppose $(\Gamma,M)$ is the following three-node E-GCM graph:}  
\hspace*{0.1in}
\parbox[c]{0.5in}{
\setlength{\unitlength}{0.75in}
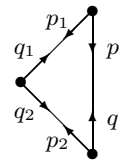
\begin{picture}(0.6,1.2)
            \put(0,0.6){\circle*{0.075}}
            \put(0.5,0.1){\circle*{0.075}}
            \put(0.5,1.1){\circle*{0.075}}
            %=======
            \put(0,0.6){\line(1,1){0.5}}
            \put(0,0.6){\line(1,-1){0.5}}
            \put(0.5,0.1){\line(0,1){1}}
            %=======
            \put(0.5,0.3){\vector(0,1){0.1}}
            \put(0.6,0.3){\footnotesize $q$}
            \put(0.5,0.9){\vector(0,-1){0.1}}
            \put(0.6,0.8){\footnotesize $p$}
            \put(0,0.6){\vector(1,1){0.2}}
            \put(0.5,1.1){\vector(-1,-1){0.2}}
            \put(0,0.6){\vector(1,-1){0.2}}
            \put(0.5,0.1){\vector(-1,1){0.2}}
            \put(-0.05,0.8){\footnotesize $q_{1}$}
            \put(0.175,1){\footnotesize $p_{1}$}
            \put(-0.05,0.35){\footnotesize $q_{2}$}
            \put(0.175,0.15){\footnotesize $p_{2}$} 
\end{picture}
} 
\hspace*{0.25in}{\sl Assume that all node pairs are odd-neighborly.   
Then $(\Gamma,M)$ is not admissible.}

{\em Proof.}  
The proof is somewhat tedious. With amplitudes as depicted in 
the proposition statement, assign numbers $a$, $b$, 
and $c$ as follows:  
\hspace*{0.1in}
\parbox[c]{0.5in}{
\setlength{\unitlength}{0.75in}
\begin{picture}(0.6,1.2)
            \put(0,0.6){\circle*{0.075}}
            \put(0.5,0.1){\circle*{0.075}}
            \put(0.5,1.1){\circle*{0.075}}
            %=======
            \put(0,0.6){\line(1,1){0.5}}
            \put(0,0.6){\line(1,-1){0.5}}
            \put(0.5,0.1){\line(0,1){1}}
            %=======
            \put(0.5,0.3){\vector(0,1){0.1}}
            \put(0.5,0.9){\vector(0,-1){0.1}}
            \put(0,0.6){\vector(1,1){0.2}}
            \put(0.5,1.1){\vector(-1,-1){0.2}}
            \put(0,0.6){\vector(1,-1){0.2}}
            \put(0.5,0.1){\vector(-1,1){0.2}}
            \put(-0.15,0.58){$c$}
            \put(0.58,0.05){$b$}
            \put(0.58,1.05){$a$}
\end{picture}
} 
\hspace*{0.2in} 
Call this position $\lambda = (a,b,c)$, so $a$ is the number at node 
$\gamma_{1}$, $b$ is at node $\gamma_{2}$, and 
$c$ is at node $\gamma_{3}$.  Without loss of generality, assume that 
$pq \leq p_{1}q_{1}$ and that $pq \leq p_{2}q_{2}$.  Set 
\[
\kappa_{1} := \frac{p p_{2} + p_{1}\sqrt{pq}}{\sqrt{pq}(2-\sqrt{pq})} 
\hspace*{0.25in}\mbox{and}\hspace*{0.25in}
\kappa_{2} := \frac{q p_{1} + p_{2}\sqrt{pq}}{\sqrt{pq}(2-\sqrt{pq})}.
\]
Assume that 
$a \geq 0$, $b \geq 0$, $c \leq 0$, 
and that 
$\displaystyle \left(\kappa_{1}-\frac{p}{q_{2}\sqrt{pq}}\right)a + 
\left(\kappa_{2}-\frac{q}{q_{1}\sqrt{pq}}\right)b + c > 0$. These 
hypotheses will be referred to as 
condition ({\tt *}).  Notice that $a$ and $b$ cannot both be zero 
under condition ({\tt *}).  A justification of the following claim 
will be given at the end of the proof: 

\noindent 
{\bf Claim:}\ \ {\sl 
Under condition} ({\tt *}) 
{\sl there is a sequence of legal node firings from initial 
position $\lambda = (a,b,c)$ which results 
in the position $\lambda' = (a', b', c') = 
(\frac{-q}{\sqrt{pq}}b, \frac{-p}{\sqrt{pq}}a, \kappa_{1}a + 
\kappa_{2}b + c)$.}  

In this case, observe that $a' \leq 0$, $b' \leq 0$, and $c' > 
0$.  Now fire at node $\gamma_{3}$ to obtain the position $\lambda^{(1)} = 
(a_{1}, b_{1}, c_{1})$ with $a_{1} = 
q_{1}[\kappa_{1}a + 
(\kappa_{2}-\frac{q}{q_{1}\sqrt{pq}})b + c]$, 
$b_{1} = q_{2}[(\kappa_{1}-\frac{p}{q_{2}\sqrt{pq}})a + 
\kappa_{2}b + c]$, and $c_{1} = -(\kappa_{1}a + 
\kappa_{2}b + c)$. Now condition ({\tt *}) implies that $a_{1} > 0$, 
$b_{1} > 0$, and $c_{1} < 0$.  
At this point to see that $\lambda^{(1)} = (a_{1}, b_{1}, c_{1})$ itself meets 
condition ({\tt *}), we only need to show that 
$\displaystyle \left(\kappa_{1}-\frac{p}{q_{2}\sqrt{pq}}\right)a_{1} + 
\left(\kappa_{2}-\frac{q}{q_{1}\sqrt{pq}}\right)b_{1} + c_{1} > 0$.  
As a first 
step, we argue that ({\em i}) $q_{1}(\kappa_{1}-\frac{p}{q_{2}\sqrt{pq}}) \geq 
1$ and that ({\em ii}) $q_{2}(\kappa_{2}-\frac{q}{q_{1}\sqrt{pq}}) \geq 1$.  
We only show ({\em i}) since ({\em ii}) follows by similar 
reasoning.  
(From the inequalities ({\em i}) and ({\em ii}), a third inequality 
({\em iii}) follows 
immediately:  
$q_{1}(\kappa_{1}-\frac{p}{q_{2}\sqrt{pq}}) + 
q_{2}(\kappa_{2}-\frac{q}{q_{1}\sqrt{pq}}) - 1 > 0$.)  
For the first of the inequalities ({\em i}), note that since $1 \leq pq$, 
then $2-\sqrt{pq} \leq pq$.  Since $pq \leq p_{2}q_{2}$, then 
$2-\sqrt{pq} \leq p_{2}q_{2}$.  (Similarly $2-\sqrt{pq} \leq 
p_{1}q_{1}$.) Thus 
$\frac{p_{2}q_{2}}{2-\sqrt{pq}}-1 \geq 0$, and hence 
$\frac{p_{2}}{2-\sqrt{pq}} - \frac{1}{q_{2}} \geq 0$.  Therefore, 
$\frac{q_{1}pp_{2}}{\sqrt{pq}(2-\sqrt{pq})} - 
\frac{q_{1}p}{q_{2}\sqrt{pq}} \geq 0$.  Since 
$\frac{p_{1}q_{1}\sqrt{pq}}{\sqrt{pq}(2-\sqrt{pq})} \geq 1$, then 
$\frac{q_{1}pp_{2}}{\sqrt{pq}(2-\sqrt{pq})} + 
\frac{p_{1}q_{1}\sqrt{pq}}{\sqrt{pq}(2-\sqrt{pq})} - 
\frac{q_{1}p}{q_{2}\sqrt{pq}} \geq 1$. From this we get 
$q_{1}(\kappa_{1}-\frac{p}{q_{2}\sqrt{pq}}) \geq 1$, which is ({\em i}). 
The following identity is easy to verify: 

\begin{eqnarray*}
\left(\kappa_{1}-\frac{p}{q_{2}\sqrt{pq}}\right)a_{1} + 
\left(\kappa_{2}-\frac{q}{q_{1}\sqrt{pq}}\right)b_{1} + c_{1} &  & 
\hspace*{3in}
\end{eqnarray*}
\begin{eqnarray*} 
\hspace*{0.75in} & = & 
\left(\kappa_{1}-\frac{p}{q_{2}\sqrt{pq}}\right)\left[
q_{1}\left(\kappa_{1}-\frac{p}{q_{2}\sqrt{pq}}\right) + 
q_{2}\left(\kappa_{2}-\frac{q}{q_{1}\sqrt{pq}}\right) - 1\right]\,\! 
a_{1}\\
\hspace*{0.75in} &  & +\   
\left(\kappa_{2}-\frac{q}{q_{1}\sqrt{pq}}\right)\left[
q_{1}\left(\kappa_{1}-\frac{p}{q_{2}\sqrt{pq}}\right) + 
q_{2}\left(\kappa_{2}-\frac{q}{q_{1}\sqrt{pq}}\right) - 1\right]\,\! 
b_{1}\\ 
\hspace*{0.75in} &  & +\ 
\left[q_{1}\left(\kappa_{1}-\frac{p}{q_{2}\sqrt{pq}}\right) + 
q_{2}\left(\kappa_{2}-\frac{q}{q_{1}\sqrt{pq}}\right) - 1\right]\,\! 
c_{1}\\ 
\hspace*{0.75in} &  & +\  
\frac{p}{q_{2}\sqrt{pq}}\left[
q_{1}\left(\kappa_{1}-\frac{p}{q_{2}\sqrt{pq}}\right) - 1\right]\,\! 
a_{1} + 
\frac{q}{q_{1}\sqrt{pq}}\left[
q_{2}\left(\kappa_{2}-\frac{q}{q_{1}\sqrt{pq}}\right) - 1\right]\,\! 
b_{1}
\end{eqnarray*}

Now the inequalities ({\em i}), ({\em ii}), and ({\em iii}) 
of the previous paragraph together with the 
inequality $\displaystyle \left(\kappa_{1}-\frac{p}{q_{2}\sqrt{pq}}\right)a + 
\left(\kappa_{2}-\frac{q}{q_{1}\sqrt{pq}}\right)b + c > 0$ 
from condition ({\tt *}) imply that 
$\displaystyle \left(\kappa_{1}-\frac{p}{q_{2}\sqrt{pq}}\right)a_{1} + 
\left(\kappa_{2}-\frac{q}{q_{1}\sqrt{pq}}\right)b_{1} + c_{1} > 0$, as 
desired.  This means that position 
$\lambda^{(1)} = (a_{1}, b_{1}, c_{1})$ meets condition ({\tt *}) and 
none of its numbers are zero.  
In view of our {\bf Claim}, we may apply to position $\lambda^{(1)}$ a 
legal sequence of node 
firings followed by firing node $\gamma_{3}$ as before to obtain a 
position $\lambda^{(2)} = (a_{2}, b_{2}, c_{2})$ that meets 
condition ({\tt *}) with none of its numbers zero, etc.  
So from any such $\lambda = (a,b,c)$ we have a divergent game 
sequence.  In view of inequalities ({\em i}) and ({\em ii}), 
the fundamental positions 
$\omega_{1} = (1,0,0)$ and $\omega_{2} = (0,1,0)$ meet condition 
({\tt *}).  The fundamental position $\omega_{3} = (0,0,1)$ does not 
meet condition ({\tt *}). However, by firing at node $\gamma_{3}$ we 
obtain the position $(q_{1}, q_{2}, -1)$, which meets  
condition ({\tt *}) by inequality ({\em iii}).  
Thus from any fundamental position there is a 
divergent game sequence, and so by \NotMarsFriendlyLemma\ of 
\cite{DonEnumbers} the 
three-node E-GCM graph we started with is not admissible. 

To complete the proof we must justify our {\bf Claim}.   Beginning with 
position $\lambda = (a,b,c)$ under condition ({\tt *}), we propose to 
fire at nodes 
$\gamma_{1}$ and $\gamma_{2}$ in alternating order until this is no 
longer possible.  
We assert that the resulting position will be 
$\lambda' = (a', b', c') = 
(\frac{-q}{\sqrt{pq}}b, \frac{-p}{\sqrt{pq}}a, \kappa_{1}a + 
\kappa_{2}b + c)$. 
There are three cases to consider: (I), $a$ and $b$ are both 
positive, (II), $a>0$ and $b=0$, and (III), $a=0$ and $b>0$.  
For (I), 
we wish to show that $(\gamma_{1},\gamma_{2},\ldots,\gamma_{1})$ of 
length $m_{12}$ is a sequence of legal node firings.  That is, we 
must check that 
\begin{eqnarray} 
\langle (s_{2}s_{1})^{k}.\lambda , \alpha_{1} \rangle & = & 
\langle \lambda , (s_{1}s_{2})^{k}.\alpha_{1} \rangle \ > \ 0 
\hspace*{0.25in} \mbox{for}\ \ 0 \leq k \leq (m_{12}-1)/2, \mbox{ and}\\
\langle s_{1}(s_{2}s_{1})^{k}.\lambda , \alpha_{2} \rangle & = & 
\langle \lambda , s_{1}(s_{2}s_{1})^{k}.\alpha_{2} \rangle \ > \ 0 
\hspace*{0.25in} \mbox{for}\ \ 0 \leq k < (m_{12}-1)/2 
\end{eqnarray}
For (II), 
we wish to show that 
$(\gamma_{1},\gamma_{2},\ldots,\gamma_{1},\gamma_{2})$ of 
length $m_{12}-1$ is a sequence of legal node firings.  That is, we 
must check that 
\begin{eqnarray} 
\langle (s_{2}s_{1})^{k}.\lambda , \alpha_{1} \rangle & = & 
\langle \lambda , (s_{1}s_{2})^{k}.\alpha_{1} \rangle \ > \ 0 
\hspace*{0.25in} \mbox{for}\ \ 0 \leq k < (m_{12}-1)/2, \mbox{ and}\\
\langle s_{1}(s_{2}s_{1})^{k}.\lambda , \alpha_{2} \rangle & = & 
\langle \lambda , s_{1}(s_{2}s_{1})^{k}.\alpha_{2} \rangle \ > \ 0 
\hspace*{0.25in} \mbox{for}\ \ 0 \leq k < (m_{12}-1)/2 
\end{eqnarray}
For (III), 
we wish to show that 
$(\gamma_{2},\gamma_{1},\ldots,\gamma_{2},\gamma_{1})$ of 
length $m_{12}-1$ is a sequence of legal node firings.  That is, we 
must check that 
\begin{eqnarray} 
\langle (s_{1}s_{2})^{k}.\lambda , \alpha_{2} \rangle & = & 
\langle \lambda , (s_{2}s_{1})^{k}.\alpha_{2} \rangle \ > \ 0 
\hspace*{0.25in} \mbox{for}\ \ 0 \leq k < (m_{12}-1)/2, \mbox{ and}\\
\langle s_{2}(s_{1}s_{2})^{k}.\lambda , \alpha_{1} \rangle & = & 
\langle \lambda , s_{2}(s_{1}s_{2})^{k}.\alpha_{1} \rangle \ > \ 0 
\hspace*{0.25in} \mbox{for}\ \ 0 \leq k < (m_{12}-1)/2 
\end{eqnarray}

To address (1) through (6), we consider matrix representations 
for each $S_{i} := \sigma_{M}(s_{i})$ (where  
$i = 1,2,3$) under the 
representation $\sigma_{M}$.  
With respect to the ordered basis $\mathfrak{B} = 
(\alpha_{1}, \alpha_{2}, \alpha_{3})$ for $V$ we have 
$X_{1} := [S_{1}]_{\mathfrak{B}} = 
\left(\begin{array}{ccc}-1 & p & p_{1}\\ 0 & 1 & 0\\ 0 & 0 & 
1\end{array}\right)$ and 
$X_{2} := [S_{2}]_{\mathfrak{B}} = 
\left(\begin{array}{ccc}1 & 0 & 0\\ 
q & -1 & p_{2}\\
0 & 0 & 1\end{array}\right)$, and so 
\[X_{1,2} := [S_{1}S_{2}]_{\mathfrak{B}} = X_{1}X_{2} = 
\left(\begin{array}{ccc}
pq-1 & -p & p_{2}p + p_{1}\\
q & -1 & p_{2}\\
0 & 0 & 1
\end{array}\right)\]
and 
\[X_{2,1} := [S_{2}S_{1}]_{\mathfrak{B}} = X_{2}X_{1} = 
\left(\begin{array}{ccc}
-1 & p & p_{1}\\
-q & pq-1 & p_{1}q + p_{2}\\
0 & 0 & 1
\end{array}\right).\]
For (1) through (6) above, we need to understand 
$X_{1,2}^{k}${\footnotesize 
$\left(\begin{array}{c}1\\ 0\\ 0\end{array}\right)$},  
$X_{2}X_{1,2}^{k}${\footnotesize 
$\left(\begin{array}{c}1\\ 0\\ 0\end{array}\right)$},  
$X_{2,1}^{k}${\footnotesize 
$\left(\begin{array}{c}0\\ 1\\ 0\end{array}\right)$}, and  
$X_{1}X_{2,1}^{k}${\footnotesize 
$\left(\begin{array}{c}0\\ 1\\ 0\end{array}\right)$}.   
Set $\theta := \pi/m_{12}$.  
Then we can write $X_{1,2} = PDP^{-1}$ for nonsingular $P$ and 
diagonal matrix $D$ as in 
\[
\frac{1}{q(e^{2i\theta}-e^{-2i\theta})}
\left(\begin{array}{ccc}
e^{2i\theta}+1 & e^{-2i\theta}+1 & p_{2}p + 2p_{1}\\
q & q & p_{1}q+2p_{2}\\
0 & 0 & 4-pq
\end{array}\right)
\left(\begin{array}{ccc}
e^{2i\theta} & 0 & 0\\
0 & e^{-2i\theta} & 0\\
0 & 0 & 1
\end{array}\right)
\left(\begin{array}{ccc}
q & -e^{-2i\theta}-1 & C_{1}\\
-q & e^{2i\theta}+1 & C_{2}\\
0 & 0 & C_{3}
\end{array}\right)
\]
where 
$C_{1} = [-q(p_{2}p+2p_{1})+(e^{-2i\theta}+1)(p_{1}q+2p_{2})]/(4-pq)$, 
$C_{2} = [q(p_{2}p+2p_{1})-(e^{2i\theta}+1)(p_{1}q+2p_{2})]/(4-pq)$, and 
$C_{3} = q(e^{2i\theta}-e^{-2i\theta})/(4-pq)$.  With some work we can 
calculate $X_{1,2}^{k}$, which results in 
\[X_{1,2}^{k} = \left(\begin{array}{ccc}
\frac{\sin(2(k+1)\theta)+\sin(2k\theta)}{\sin(2\theta)} & 
\frac{-p\sin(2k\theta)}{\sin(2\theta)} & C_{1}'\\
\frac{q\sin(2k\theta)}{\sin(2\theta)} & 
\frac{-\sin(2k\theta)-\sin(2(k-1)\theta)}{\sin(2\theta)} & 
C_{2}'\\
0 & 0 & 1
\end{array}\right),\]
with \[C_{1}' = -\frac{p_{2}p+2p_{1}}
{4-pq}\left[\frac{\sin(2(k+1)\theta)+\sin(2k\theta)}
{\sin(2\theta)} 
- 1\right] + 
\frac{p(p_{1}q+2p_{2})\sin(2k\theta)}{(4-pq)\sin(2\theta)}\]  
and \[C_{2}' = -\frac{q(p_{2}p+2p_{1})\sin(2k\theta)}{(4-pq)\sin(2\theta)} + 
\frac{p_{1}q+2p_{2}}
{4-pq}\left[\frac{\sin(2k\theta)+\sin(2(k-1)\theta)}
{\sin(2\theta)} + 1\right].\] 
Similar reasoning (or simply interchanging the roles of $\alpha_{1}$ 
and $\alpha_{2}$ in the preceding calculations, or noting that 
$X_{2,1}^{k} = (X_{1,2}^{-1})^{k} = X_{1,2}^{-k}$\, ) shows that 
\[X_{2,1}^{k} = \left(\begin{array}{ccc}
\frac{-\sin(2k\theta)-\sin(2(k-1)\theta)}{\sin(2\theta)} & 
\frac{p\sin(2k\theta)}{\sin(2\theta)} & C_{1}''\\
\frac{-q\sin(2k\theta)}{\sin(2\theta)} & 
\frac{\sin(2(k+1)\theta)+\sin(2k\theta)}{\sin(2\theta)} & 
C_{2}''\\
0 & 0 & 1
\end{array}\right),\]
with 
\[C_{1}'' =   
\frac{p_{2}p+2p_{1}}
{4-pq}\left[\frac{\sin(2k\theta)+\sin(2(k-1)\theta)}
{\sin(2\theta)} + 1\right] - 
\frac{p(p_{1}q+2p_{2})\sin(2k\theta)}{(4-pq)\sin(2\theta)}\] and 
\[C_{2}'' = \frac{q(p_{2}p+2p_{1})\sin(2k\theta)}{(4-pq)\sin(2\theta)} - 
\frac{p_{1}q+2p_{2}}{4-pq}\left[\frac{\sin(2(k+1)\theta)+\sin(2k\theta)}
{\sin(2\theta)} 
- 1\right].\]  
Then 
\[X_{2}X_{1,2}^{k} = \left(\begin{array}{ccc}
\frac{\sin(2(k+1)\theta)+\sin(2k\theta)}{\sin(2\theta)} & 
\frac{-p\sin(2k\theta)}{\sin(2\theta)} & C_{1}'\\
\frac{q\sin(2(k+1)\theta)}{\sin(2\theta)} &
\frac{(1-pq)\sin(2k\theta)+\sin(2(k-1)\theta)}{\sin(2\theta)} & 
qC_{1}' - C_{2}' + p_{2}\\
0 & 0 & 1
\end{array}\right)\]
and 
\[X_{1}X_{2,1}^{k} = \left(\begin{array}{ccc}
\frac{(1-pq)\sin(2k\theta)+\sin(2(k-1)\theta)}{\sin(2\theta)} & 
\frac{p\sin(2(k+1)\theta)}{\sin(2\theta)} & -C_{1}'' + pC_{2}'' + p_{1}\\
\frac{-q\sin(2k\theta)}{\sin(2\theta)} & 
\frac{\sin(2(k+1)\theta)+\sin(2k\theta)}{\sin(2\theta)} & 
C_{2}''\\
0 & 0 & 1
\end{array}\right).\] 
Now we can justify (1) through (6).  For example, for (4) we see that 
since $X_{1}X_{2,1}^{k}${\footnotesize 
$\left(\begin{array}{c}0\\ 1\\ 0\end{array}\right)$} is the second 
column of the matrix $X_{1}X_{2,1}^{k}$, 
then $\langle \lambda , s_{1}(s_{2}s_{1})^{k}.\alpha_{2} \rangle = 
a\frac{p\sin(2(k+1)\theta)}{\sin(2\theta)}$, which is positive 
since $a > 0$, $p > 0$, and (recalling that $m_{12}$ is odd) 
$2(k+1) < m_{12}$.  

Then the proposed firing sequence for each of cases (I), (II), and 
(III) is legal.  To see in case (I) 
that the resulting position is the claimed 
$\lambda' = (a', b', c') = 
(\frac{-q}{\sqrt{pq}}b, \frac{-p}{\sqrt{pq}}a, \kappa_{1}a + 
\kappa_{2}b + c)$, 
we need to calculate 
$\langle s_{1}(s_{2}s_{1})^{k}.\lambda , \alpha_{i} \rangle = 
\langle \lambda , s_{1}(s_{2}s_{1})^{k}.\alpha_{i} \rangle$ for each 
of $i = 1,2,3$, where $k$ is now $(m_{12}-1)/2$.  With patience one 
can confirm that 
\[X_{1}X_{2,1}^{k} = 
\left(\begin{array}{ccc}
0 & 
-p/\sqrt{pq} & \kappa_{1}\\
-q/\sqrt{pq} & 0 &  
\kappa_{2}\\
0 & 0 & 1
\end{array}\right),\] from which the claim follows.  Similar 
computations confirm the claim for cases (II) and (III).\hfill\QED

\vspace*{-0.15in}

%=============================================
% Bibliography
%=============================================
\renewcommand{\baselinestretch}{1}

\end{document}